\documentclass{amsart}
\usepackage[]{amsmath, amsthm, amsfonts, verbatim, amssymb}
\numberwithin{equation}{section}

\def\cM{{\mathcal M}}

\def\cT{{\mathcal T}}

\def\bR{{\mathbb R}}

\def\bC{{\mathbb C}}

\newtheorem {theo}{Theorem}
\newtheorem {coro}{Corollary}
\newtheorem {lemm}{Lemma}

\newtheorem {defi}{Definition}

\newtheorem {prop}{Proposition}

\def\tM{\tilde M}

\def\tell{\widetilde \ell}
\def\bs{\bigskip}
\def\ms{\medskip}
\def\ni{\noindent}
\def\td{\nabla}
\def\pd{\partial}

\def\oW{\overline W}

\def\oI{\overline I}

\def\oj{\overline j}
\def\oi{\overline i}

\def\ou{\overline u}
\def\ov{\overline v}
\def\ow{\overline w}

\def\oM{\overline M} 
\def\oA{\overline A}

\def\oIq1{\oI_1\cdots\oI_{q-1}}
\def\oIq2{\oI_1\cdots\oI_{q-2}}
\def\op{\overline \partial}

\def\vol{\mbox{vol}}

\def\ddz{\frac{\pd}{\pd z}}
\def\ddzi{\frac{\pd}{\pd z^i}}

\def\ddzk{\frac{\pd}{\pd z^k}}
\def\oddzk{\frac{\pd}{\pd \bar z^k}}

\def\ddbwj{\frac{\pd}{\op w^j}}
\def\ddw{\frac{\pd}{\pd w}}
\def\ddz1{\frac{\pd}{\pd z^1}}
\def\ddw1{\frac{\pd}{\pd w^1}}
\def\oddw1{\overline{\frac{\pd}{\pd w_1}}}

\def\oddzj{\overline{\frac{\pd}{\pd z^j}}}
\def\oddwi{\overline{\frac{\pd}{\pd w^i}}}
\def\oddwj{\overline{\frac{\pd}{\pd w^j}}}
\def\ddwi{\frac{\pd}{\pd w^i}}

\def\td{\nabla}
\begin{document}

\title[Geometry of domains with the uniform squeezing property]
{Geometry of domains with the uniform squeezing property}
\author[Sai-Kee Yeung]
{Sai-Kee Yeung}

\address[]{Mathematics Department, Purdue University, West Lafayette, IN  47907
USA} \email{yeung@math.purdue.edu}

\thanks{\noindent{\it \noindent Key Words: bounded domains, pseudoconvexity,
K\"ahler-Einstein
metric, K\"ahler hyperbolic,  moduli space\\
1991 Mathematics Subject Classification.}  Primary  32G15, 53C55, 55N99\\The author was partially supported by grants from the National Science Foundation
and the National Security Agency}

\maketitle

\bigskip
\noindent
{\bf Abstract} {\it We introduce the notion of domains with uniform
squeezing property, study various analytic and
geometric properties of such domains and show that they cover many
interesting examples, including Teichm\"uller spaces and Hermitian symmetric
spaces of non-compact type.  The properties
supported by such manifolds include pseudoconvexity, hyperconvexity, 
K\"ahler-hyperbolicity, vanishing of cohomology groups and quasi-isometry
of various invariant metrics.  It also leads to nice geometric properties 
for manifolds covered
by bounded domains and a simple criterion to provide positive examples
to a problem of Serre about Stein properties of holomorphic fiber bundles.}

\bs
\begin{center}
{\bf \S0 Introduction}
\end{center}

The purpose of this article is to introduce a class of bounded domains
in $\bC^n$ which on one hand is sufficiently general to 
include interesting classes of examples and
on the other hand leads to interesting 
analytic and geometric properties.  

\begin{defi}
Denote by $B_r(x)$ a ball of radius $r$ in $\bC^n.$
Let $0<a<b<\infty$ be positive constants.
A bounded domain $M$ in $\bC^n$ for some $n>1$ 
is said to have the uniform squeezing property,
or more precisely, $(a,b)$-uniform squeezing property
if there exist constants $a$ and $b$,
such that for each point $x\in M,$ there exists an 
embedding $\varphi_x:M\rightarrow \bC^n$ with $\varphi_x(x)=0$ and 
$B_a(\varphi_x(x))\subset\varphi_x(M)\subset B_b(\varphi_x(x)).$
We call the corresponding coordinate system a uniform squeezing
coordinate system or, more precisely, $(a,b)$-uniform coordinate system.
\end{defi}

Even though the definition is very simple and appears to be rather
restrictive, it in fact includes lots of interesting examples.

\begin{prop}
Examples of bounded domains with the uniform squeezing property include
the followings,\\
(a). bounded homogeneous domains,\\
(b).  bounded strongly convex domains,\\
(c). bounded domains which cover a compact K\"ahler manifold, and\\
(d). Teichm\"uller spaces $\cT_{g,n}$ of hyperbolic Riemann surfaces
of genus $g$ with $n$ punctures.
\end{prop}

We remark that Hermitian symmetric spaces of non-compact type constitute
an important subclass of both (a) and (c).  The former follows from
the Harish-Chandra realization of such symmetric spaces, and the latter follows from
existence of cocompact arithemtic lattices associated to the automorphism groups of
the symmetric spaces (cf. [B]).

Our main objective is to show that domains with uniform
squeezing properties support many interesting geometric and analytic 
properties.
The first observation is about pseudoconvexity of such domains.

\begin{theo} Let $M$ be a domain with the uniform squeezing property.
Then the following conclusions are valid.\\
(a). The Bergman metric of $M$ is complete.\\
(b). $M$ is a pseudoconvex domain.\\
(c). There exists a complete K\"ahler-Einstein metric on $M$.
\end{theo}

The second observation is about the behavior of invariant metrics on
such domains.  On a general bounded domain, there are three well-known intrinsic metrics
which are invariant under a biholomorphism, namely, the 
Kobayashi metric, the Carath\'eodory metric and
the Bergman metric.  There is a fourth one when the bounded domain is pseudoconvex, viz., the K\"ahler-Einstein
metric.  
We denote the metrics by
$g_K, g_C, g_B$ and $g_{KE}$ respectively.  For a K\"ahler metric $g$ on
$M$, we denote by $R^g$ its curvature tensor and $\td^g$ the Riemannian
connection.  For tangent vectors $X_1,\cdots, X_N,$ we denote $\td_{X_1}^g\cdots\td_{X_N}^g$
by $\td_{X_1,\cdots 
X_N}^{g}.$  Furthermore, $\td^{g}_{i_1,\cdots i_N}$ denotes the covariant derivatives with
respect to the coordinate vectors.
We also normalized the K\"ahler-Einstein metric so that
$Ric(g_{KE})=-2(n+1).$

\begin{theo}
Let $M$ be a bounded domain with
$(a,b)$-uniform squeezing property.\\
(a).  The invariant metrics
 $g_K, g_C, g_B$ and $g_{KE}$ are quasi-isometric.  Furthermore
\begin{eqnarray*}
\frac ab g_K\leqslant&g_C&\leqslant g_K\\
\frac ab g_K\leqslant&g_B&\leqslant [\frac{2\pi}{a^3}(\frac {2b}a)^n]^2g_K\\
\frac{a^2}{b^2n}g_K\leqslant&g_{KE}&\leqslant(\frac {b^{4n-2}n^{n-1}}a)g_K
\end{eqnarray*}
\ni(b). 
There exist constants $c^{g_B}_n$ and $c^{g_{KE}}_n$
such that $\Vert \td^{g_B}_{i_1,\cdots i_N}R_{g_B}\Vert_{g_B}\leqslant c^{g_B}_N$
and $\Vert \td^{g_{KE}}_{i_1,\cdots i_N}R_{g_{KE}}\Vert_{g_{KE}}\leqslant c^{g_{KE}}_N$ for any covariant derivatives $\td^{g_{B}}_{i_1,\cdots i_N}$ 
and $\td^{g_{KE}}_{i_1,\cdots i_n}$of $g_B$ and $g_{KE}$ respectively.\\
(c).  Let $X_1,\cdots, X_N$ be $N$ tangent vectors of unit length with respect to a metric $g_1$ at $x\in M.$
Then  $\Vert\td_{X_1,\cdots 
X_N}^{g_1}g_B-\td_{X_1,\cdots X_N}^{g_1}g_{KE}\Vert_{g_1}\leqslant c,$
for some constant $c$ depending on $N,$ where $g_1=g_B$ or $g_{KE}.$
\\
(d).  Both of $g_{B}$ and $g_{KE}$ are geometrically finite in the sense 
that
they are complete with bounded curvature and the injectivity radius is
 bounded from below uniformly on $M$.\\
(e). Both of $g_{B}$ and $g_{KE}$ are K\"ahler-hyperbolic.\\
(f). $M$ is hyperconvex.
\end{theo}

Here we recall that a K\"ahler manifold $(X,\omega)$
is K\"ahler-hyperbolic if on its
universal covering, $\omega$ can be written as $dh$ where $h$ is bounded
uniformly when measured with respect to $\omega.$ 
$M$ is hyperconvex if there exists a plurisubharmonic exhaustion
function bounded from above on $M.$

The followings are some well-known consequences of K\"ahler-hyperbolicity
in Theorem 1 and Theorem 2.

\begin{coro} Let $M$ be a uniformly squeezed manifold.
Let $g$ be $g_B$ or $g_{KE}.$\\  
(a) The reduced $L^2$-cohomology groups of $M$ with respect to $g=g_B$ 
or $g_{KE}$ satisfies
$\dim(H^0_{(2)}(M))=\infty$
 and $\dim(H^i_{(2)}(M))=0$ for all $i>0.$\\
(b) The first eigenvalues of the Beltrami Laplacian operators
 $\Delta_{g_B}$ and $\Delta_{g_{KE}}$ on smooth functions on $M$ with 
respect to $g_B$ and $g_{KE}$ respectively 
are both bounded from below by $0.$\\
(c) The volume with respect to either $g=g_B$ or
$g_{KE}$
of any relatively compact complex submanifold
with boundary 
$N\subset M$
of complex dimension $k$ satisfies $\vol_{2k}(N)\leqslant C\cdot \vol_{2k-1}(\pd N)$
for some  constant $C>0.$
\end{coro}

We define a lattice on $M$ to be a discrete group acting properly discontinuously as biholomorphisms
on $M.$ 

\begin{coro}  Assume that $\Gamma$ is a torsion-free lattice on $M$ which admits a
uniform squeezing coordinate system.  Then a compact quotient
 $N=M/\Gamma$ has to be a projective
algebraic variety of general type.
A non-compact quotient $N=M/\Gamma$ which has 
finite volume with respect to the invariant Bergman metric has
to be a quasi-projective
variety of log-general type.
\end{coro}

Another direct consequence of Theorem 2b and a result of Mok-Yau in [MY] is the following estimates on the growth of Bergman
kernel.

\begin{coro} 
Let $M$ be a bounded domain with
the uniform squeezing property.  Denote by $d=d(z,\pd\Omega)$ the Euclidean distance of $z\in \Omega$
from the boundary $\pd \Omega$ of $\Omega.$  Then
$K(z,z)\geqslant \frac c{d^2(-\log d)^2}$ for some constant $c>0.$
\end{coro}

Let us now focus on the 
applications of the above results to more specific compact or non-compact
 manifolds.

\begin{theo} 
Assume that $N$ is a compact complex manifold of complex dimension $n$
whose universal covering is a bounded domain in $\bC^n.$  Then the following 
properties hold. \\
(a). $N$ is projective algebraic.\\
(b). There exists a K\"ahler-Einstein metric on $N$.\\
(c). $(-1)^n\chi(N)>0.$\\
(d). $H^0(N,2K)$ is non-trivial, where $K$ is the canonical line bundle on
$N$.\\
(e). The universal covering of $N$ is Stein.
\end{theo}

The result can be considered as a support for a conjecture of Shafarevich,
which states that the universal covering of a projective algebraic  
variety is holomorphically convex (cf. [Ko]).  The assumption is stronger but the
projective algebraicity is obtained for free.  On the other hand, it also
shows that if we try to formulate a conjecture for the uniformization
of a compact complex manifold by a bounded domain, it should include 
topological and analytic conditions such as those stated in (c) and (d).
Properties in (c) is along the line of conjectures of Hopf, Chern and Singer in Riemannian 
geometry,
a consequence of those is that the Euler characteristic
of a compact Riemannian manifold of even dimension $2n$ with non-positive Riemannian sectional
curvature satisfies (c) (cf. [Gr]).  Note that a compact torus is flat and its Euler characteristic is equal to zero.

\bs 
As a consequence of a result of Stehl\'e, Theorem 2f provides the following 
simple criterion for
positive solutions to a problem of Serre [Se], who asked whether a holomorphic
fiber bundle with Stein base and Stein fibers are Stein.

\begin{coro}
Suppose $\pi:T\rightarrow B$ is a locally trivial holomorphic fiber space
for which the base $B$ is a Stein space and the fibers satisfy the uniform
squeezing properties.  Then $T$ is also Stein.
\end{coro}

\bs
As an application of Theorem 1 and 2 to non-compact complex manifolds of
finite volume with respect to some invariant metric, we consider
moduli space of possibly punctured curves as an example.

\begin{theo}
Let $g,n\geqslant 0$ and $2g-2+n>0,$ so that the complement of $n$ punctures
of a compact Riemann surface of genus $g$ gives a hyperbolic Riemann surface.
Let $\cM_{g,n}$ be the moduli space of such hyperbolic Riemann surfaces.  Let $\cT_{g,n}$ be the
corresponding Teichm\"uller space.\\  
(a). $g_K, g_C, g_B$ and $g_{KE}$ are quasi-isometric on $\cM_{g,n}$\\
(b). 
For $g_B$ and $g_{KE}$, any order of 
 covariant derivative of the curvature tensor of the metric is uniformly
bounded on $\cM_{g,n}.$  As a consequence, for any set of unit vectors $\{ X_1,\cdots 
X_N\}$ measured with respect to $g_1,$ the difference $\Vert\td_{X_1,\cdots 
X_N}^{g_1}R^{g_1}-\td_{X_1,\cdots X_N}^{g_2}R^{g_2}\Vert_{g_1}$ is bounded
for any $g_1, g_2$ chosen among $g_B,$ $g_{KE}$, where $\td_{X_1,\cdots X_N}^{g}$ denotes the
covariant derivatives of a metric $g$ with respect to the vectors $X_1,\cdots, X_n.$ \\
(c). The Teichm\"uller space $\cT_{g,n}$ is K\"ahler-hyperbolic with
respect to both $g_B$ and $g_{KE}.$\\
(d). $\cT_{g,n}$ is hyper-convex.\\
(f). $\cM_{g,n}$ is quasi-projective of log-general type and 
the Euler-Poincar\'e characteristic satisfies 
$(-1)^n\chi(\cM_{g,n})>0.$

\end{theo}

Except for the statements related to the K\"ahler-hyperbolicity of the K\"ahler-Einstein
metric and estimates on higher order quasi-isometry of the metrics
$g_B$ and $g_{KE}$, most of the results in Theorem 4 can be obtained for example
by combining results in [Y3]
and [Y4], but the proofs there rely on many well-known and diverse results.
In this paper, all these properties
except quasi-projectivity of moduli space of curves  
are derived solely from the existence of uniform
squeezing coordinates, which is provided classically
by the Bers Embedding (cf. [Ga]).

\bs
Overall we remark that parts of the results in this paper have been obtained 
for some specific examples mentioned in Proposition 1.
In particular, K\"ahler-hyperbolicity of locally Hermitian
symmetric spaces with respect to
the Bergman metric is explained in [Gr], of bounded homogeneous space with
respect to the Bergman metric is
proved in [Do], of moduli space of curves with respect to a metric
constructed by McMullen
is proved by McMullen in [Mc], where the metric is also shown to have
many nice properties such as geometric finiteness and quasi-isometry
to the Kobayashi metric.
The vanishing of the cohomology groups $h^i_{(2)}, i<n,$ for  K\"ahler-hyperbolic manifolds was proved
by Gromov in [Gr].

There is a vast amount of literature related to Theorem 4.
A precise formula for the Euler characteristic of the moduli space of curves
was given by Harer-Zagier.  $\cM_{g,n}$ was first shown to be pseudoconvex
or a domain of holomorphy through the work of Bers and Ehrenpreis. 
Hyperconvexity of $\cM_{g,n}$ had been
proved by Krushkal, and more recently in [Y2] by a different and 
more geometric 
 way 
using results of Wolpert.  It follows from classical results of
Baily, Deligne-Mumford and Knudson-Mumford that a moduli space of
curves is quasi-projective. 
Quasi-isometry among invariant
metrics on $\cM_{g,n}$ were obtained through the contributions of Chen,
Liu-Sun-Yau and Yeung.
We refer the readers to [W] and [Y4] for more details.

\bs
In a sequel to the present paper, we will explain how the set-up 
of the paper can be used to prove subelliptic estimates for solutions
to the $\op$ equations
on uniform squeezing domains, which include Teichm\"uller spaces, where
the problem remained open in the past due to a lack of
description of the boundary of Bers Embedding.

\bs
Parts of the work  were finished while the author was visiting
the Korea Institute for Advanced Study, the Osaka University, and the University of Hong Kong. The author would like
to thank Jun-Muk Hwang, Toshiki Mabuchi and Ngaiming Mok for their hospitality.
The author would also like to express his gratitude to the referee for very helpful
comments on the paper.

\bs
\begin{center}
{\bf \S1 Terminology and preliminaries }
\end{center}

Recall the following standard notations about various convexity of a domain.
A domain in $\bC^n$ is pseudoconvex if there exists a plurisubharmonic
exhaustion.  A bounded domain $\Omega=\{z\in\bC^n|r(z)<0\}$ for some 
$C^2$ function
$r(z)$ in $z=(z_1,\cdots,z_n)$ is strongly pseudoconvex if the Levi form
$\sqrt{-1}\pd\op r>0$ in a neighborhood of $\pd\Omega.$   A domain $\Omega$
in $\bC^n$ is hyperconvex if there exists a bounded plurisubharmonic 
exhaustion function.

A K\"ahler metric $\omega$ on a complex manifold $M$ is said to be 
K\"ahler-hyperbolic if on the universal covering $\tM$ of $M,$ 
the pull back of $\omega$ can be expressed as $d\eta$ for some $1$-form $\eta$ which
is bounded uniformly on $\tM$ with respect to $\omega.$

We say that two metrics $g_1$ and $g_2$ are quasi-isometric, denoted by
$g_1\sim g_2,$ if there exists a positive constant $c$ such that 
$\frac1cg_1(v,\ov)\leqslant g_2(v,\ov)\leqslant cg_1(v,\ov)$ for all
holomorphic tangent vectors $v.$

Let us now recall the various notions of invariant metrics on a general
complex manifold.

For a unit tangent vector $v\in T_xM$
 on a complex
manifold $M$, the Kobayashi and Carath\'eodory semi-metrics are defined respectively
as 
complex Finsler metrics by 
 
\begin{eqnarray*}
\sqrt{g_K(x,v)}&=&\inf\{\frac1r\vert\exists f:B_r^1\rightarrow M \mbox{ holomorphic}, 
f(0)=x, f'(0)=v\}.\\
\sqrt{g_C(x,v)}&=&\sup\{\frac1r\vert\exists h:M\rightarrow B_r^1\mbox{ holomorphic},
h(x)=0, |dh(v)|=1\},
\end{eqnarray*}
where we use $B^n_r=B^n_r(0)$ to denote a ball of radius $r$ centered at
$0$ in $\bC^n.$
Since we are considering only bounded domains in $\bC^n,$
 both $g_K$ and $g_C$ are non-degenerate
 complex Finsler metrics.  

Consider now K\"ahler-Einstein metric of constant negative scalar
curvature.  We normalize the curvature so that 
$g_{KE}$ satisfies $Ric(g_{KE})=-2(n+1)$, 
where $\omega_{KE}$ is the K\"ahler form 
associated to $g_{KE}.$  The normalization is chosen so that it
agrees with the one for the hyperbolic metric on $B^n_\bC$ of
constant holomorphic sectional curvature $-4.$   

The Bergman pseudometric $g_B$ on a general complex manifold $M$ of 
complex dimension
$N$ is a K\"ahler pseudometric with local potential given by the
coefficients of the
Bergman kernel $K(x,x).$  It is clearly 
non-degenerate for $\cT_g$.  
$g_B$ can be interpreted in the following way.

Let $f$ be a $L^2$-holomorphic $N$-form on $M$, where $\dim_\bC M=N.$
  In terms of local
coordinates $(z_1,\cdots,z_N)$ on a coordinate chart $U$, 
let $e_{K_M}=dz^1\wedge\cdots\wedge dz^N$
be a local basis of the canonical line bundle $K_M$ on $U$.
We can write $f$ as $f_Ue_{K_M}$ on $U$.
Let $f_i, i\in N$ be an 
orthonormal
basis of $L^2$-sections in $H^0_{(2)}(M,K_M).$ 
Note that from conformality, the choice is independent of the metric on $M$.
The Bergman kernel is given by 
$K(x,x)=\sum_if_i\wedge\overline{f_i}.$ 
 Let $K_U(x,x)=\sum_if_{U,i}\overline{f_{U,i}}$ be the coefficient of $K(x,x)$
in terms of the local coordinates.
The Bergman metric is given
by a K\"ahler form 
$$\omega_B=\sqrt{-1}\pd\op\log K_U(x,x)=\sqrt{-1}\frac1{K_U(x,x)^2}\sum_{i<j}(f_i\pd f_j-f_j\pd f_i)\wedge\overline{(f_i\pd f_j-f_j\pd f_i)}
,$$
which is clearly independent of the choice of a basis and $U$.
As the Bergman kernel is independent of basis, for each fixed point $x\in M,$
$$K_U(x,x)=\sup_{f\in H^0_{(2)}(M,K_M),\Vert f\Vert=1}|f_U(x)|^2,$$
where $\Vert\cdot\Vert$ stands for the $L^2$-norm.
We may assume that $\sup_{f\in H^0_{(2)}(M,K_M)}|f_U(x)|$ is realized by 
$f_x\in H^0_{(2)}(M,K)$ with $\Vert f_x\Vert=1$ so that $K_U(x,x)=|f_{x,U}(x)|^2.$
Using the fact that the Bergman kernel is independent of the 
choice of a basis again
and letting $V\in T_xM,$
$$\omega_B(V,\bar V)=
\frac1{|f_{x,U}(x)|^2}\sup_{f\in H^0_{(2)}(M,K_M),\Vert f\Vert=1, f(x)=0}|V(f_U)|^2.$$ 
Consider in particular $V=\ddzi.$  
We may also assume that the supremum for 
$|\ddzi(f_U)|^2$ among all $f\in H^0_{(2)}(M,K_M),\Vert f\Vert=1, f(x)=0$
is achieved by 
$g_{i,x}\in H^0_{(2)}(M,K_M)$
of $L^2$-norm $1.$  Hence
$\sup_{f\in H^0_{(2)}(M,K_M),\Vert f\Vert=1}|\ddzi f_U|^2=|\ddzi g_{i,x,U}(x)|^2.$
To simplify our notation, we may simply write
$$\omega_B(\ddzi,\overline{\ddzi})=
\frac1{|f_{x}(x)|^2}\sup_{f\in H^0_{(2)}(M,K_M),\Vert f\Vert=1, f(0)=0}|\ddzi(f)|^2
=\frac{|\ddzi g_{i,x}(x)|^2}{|f_{x}(x)|^2},
$$  
since the expression is clearly independent of the choice of $U$ and metric on $e_U$.

Finally let us include here two regularity estimates required for later
calculations for the convenience of the readers.
We denote by $W_{k,p}$ and 
$C_{k,\alpha}$ the spaces of functions on $B_a$ which are
bounded with respect to the Sobolev norm $\Vert\cdot\Vert_{k,p}$ and H\"older norm $|\cdot|_{k,\alpha}$ on $B_a(x)$ 
respectively.  We refer the readers to [GT] for standard notations. 

\begin{prop}(cf.[GT], page 235, 90)
Let $\Omega'\subset\subset\Omega$ be bounded domains in $\bR^{n}$ with $C^\infty$ boundary.
Let $L$ be a second order linear differential operator defined by
$$Lu=a^{ij}(x)D_{ij}u+b^i(x)D_iu+c(x)u$$
with sums over repeated indices. Let $u$ be a strong solution of the equation $Lu=f.$
\ni (a) (Calderon-Zygmund estimates) Suppose that $u\in H_1(\Omega)$ is a strong solution of $Lu=f$
with $f\in L^2(\Omega).$  
Assume that for $a^{ij},b^i,c
\in C^{0}(\bar\Omega),$
\begin{eqnarray*}
 a^{ij}v_iv_j&\ge&\lambda|v|^2\ \ \ \forall v\in\bR^{n},\\ 
|a^{ij}|, |b^i|, |c|&\leqslant& \Lambda.
\end{eqnarray*}   Then

$$\Vert u\Vert_{2,p,\Omega'}\leqslant C_1(\Vert u\Vert_{0,p,\Omega}+\Vert f\Vert_{0,p,\Omega})$$
with constant $C_1$ depending on $n,p,\lambda,\Lambda,\Omega',\Omega$ and the moduli of
continuity of $a^{ij}$ on $\Omega'.$\\
(b) (Schauder estimates)  Suppose $f\in C^\alpha(\bar\Omega)$ and 
$$|a^{ij}|_{\alpha,\Omega}, |b^i|_{\alpha,\Omega},|c|_{\alpha,\Omega}\leqslant \Lambda.$$
Then
$$|u|_{2,\alpha,\Omega'}\leqslant C_2(|u|_{0,\alpha}+|f|_{0,\alpha,\Omega}),$$
with constant $C_2$ depending on $n,\alpha,\lambda,\Lambda,\Omega$ and $\Omega'.$
\end{prop}

In our application, we will always assume that $\Omega=B_a(0)$ and $\Omega'=B_{a/2}(0)$ for a fixed value $a$, after identifying an arbitrary point
on the manifold to the origin with respect to a uniform squeezing coordinate system.  We are interested in the estimates
of the bounds and use it to show uniform bound over the manifold of
our interest instead of regularity, which is already known for general elliptic equations.
\bs

\begin{center}
{\bf \S2 Pseudoconvexity and related properties }
\end{center}

Throughout this section and \S3, \S4, we let $M$ be a bounded domain with uniform squeezing coordinates.  

\begin{lemm}
The Bergman metric $g_B$ on $M$ is a well-defined complete K\"ahler metric.
Furthermore, $g_B$ is quasi-isometric to $g_K$ as a complex Finsler metric.
\end{lemm}

\ni{\bf Proof}  The $(1,1)$-form $\omega_B$ defined in \S1 is only semi-definite in general.  
We need to show that it is in fact positive definite and gives rise to a complete metric
in our situation.

The K\"ahler form $\omega$ of the Bergman metric is given by
$$\omega_B(\ddzi,\overline{\ddzi})=\frac{|\ddzi g_{i,x}(x)|^2}{|f_{x}(x)|^2}
,$$
where $f_x$ is a function with $L^2$-norm $\Vert t\Vert=1$ realizing the supremum
of $|f(x)|$ among $L^2$-holomorphic functions 
$f\in H^0_{(2)}(M),\Vert f\Vert=1$
on $\cT,$ and $g_{i,x}$ is a holomorphic function realizing supremum of
$|\ddzi(f)|^2$ among all $f\in H^0_{(2)}(M),\Vert f\Vert=1, f(x)=0.$

From assumption  $B_{a}^n(x)\subset M\subset B_{b}^n(x),$ where $B_r^{n}(x)$
denotes a complex ball of radius $r$ centered at $x$ identified with $0$ in $\bC^n.$ 
Let $\vol_o$ denote the Euclidean volume on $\bC^n.$ 
Clearly from the Mean Value Inequality
$$(f_x(x))^2\leqslant \frac{\int_{B_{a}^{n}(x)}|f_x|^2}{\vol_o(B_{r}^{n}(x))}
\leqslant \frac{\int_{M}|f_x|^2}{\vol_o(B_{a}^{n}(x))}=\frac 1{a^{2n}\vol_o(B_{1}^{n})}.$$
The constant function 
$h_1(x)=1$ satisfies $h_1(1)=1$ and 
$$\Vert h_1\Vert^2=\vol_o(M)\leqslant \vol_o(B_{b}^{n})=b^{2n}\vol_o(B_{1}^{n}).$$

Hence $|f_x(x)|\ge\frac 1{[(b)^{2n}\vol_o(B_1^{n})]^{\frac12}}.$  We conclude that
$$[\frac {1}{a^{2n}\vol_o(B_{1}^{n})}]^{\frac12}
\ge |f_x(x)|\ge [\frac {1}{b^{2n}\vol_o(B_{1}^{n})}]^{\frac12}.$$

Let $V_i$ be the complex line generated by $\ddzi$ in $\bC^n.$
Then from Generalized Cauchy Inequality and Mean Value Inequality,
\begin{eqnarray*}
|\ddzi g_{i,x}(x)|&\leqslant&\frac{\int_{\pd(B_{\frac a2}^n(x))\cap V_i}|g_{i,x}(y)| dy}{(\frac a2)^2}
\leqslant\frac{[\int_{\pd(B_{\frac a2}^n(x))\cap V_i}|g_{i,x}(y)|^2 dy]^
{\frac12}
[2\pi\frac a2]^{\frac12}}{(\frac a2)^2} \\
&\leqslant&\frac{[\int_{\pd(B_{\frac a2}^n(x))\cap V_i}dy\int_{B_{\frac a2}^n(y)}|g_{i,x}(w)|^2d\vol_o(w)]^{\frac a2}[\frac\pi2]^{\frac 12}}
{(\frac a2)^2[\vol_o(B_{\frac a2}^n)]^{\frac12}}\\
&\leqslant&\frac{[\int_{\pd(B_{\frac a2}^n(x))\cap V_i}dy\int_{M}|g_{i,x}(w)|^2d\vol_o(w)]^{\frac12}[\frac\pi2]^{\frac12}}
{(\frac a2)^2[\vol_o(B_{\frac a2}^n)]^{\frac12}}
\leqslant\frac{\frac{\pi}2}{(\frac a2)^{2+n}[\vol_o(B_{1}^n)]^{\frac12}}.
\end{eqnarray*}

On the other hand the function $h_{i,x}=z_i$ satisfies $\ddzi h_{i,x}=1$
and $h_{i,x}(0)=0.$ 
As $\int_{B_1^n}|z_i|^2=\frac1{n+1}\vol(B_1^n)$, we know
that 
$$\Vert h_{i,x}\Vert^2\leqslant\int_{B_{b}^n}|z_i|^2
=b^{2n+1}\int_{B_1^n}|z_i|^2\leqslant \frac1{n+1}b^{2n+1}\vol_o(B_{1}^{n}).$$
Hence the function $k_{i,x}:=\frac{h_{i,x}}{\Vert h_{i,x}\Vert}$ satisfies
$|\ddzi k_{i,x}|=\frac {\sqrt{n+1}}{[b^{2n+1}\vol_o(B_{1}^{n})]^{\frac12}}$, $k_{i,x}(0)=0$ and 
$\Vert k_{i,x}\Vert^2=1.$

We conclude as before that
$$\frac{\frac{\pi}2}{(\frac a2)^{n+2}[\vol_o(B_{1}^n)]^{\frac12}}\ge |\ddzi g_{i,x}(x)|
\ge\frac{\sqrt{n+1}}{b^{n+\frac12}[\vol_o(B_{1}^{n})]^{\frac12}}.$$ 

Combining the above estimates for $f_x(x)$ and $g_{i,x}(x),$ we arrive at
$$\frac{2\pi}{a^2}(\frac {2b}a)^n\ge \sqrt{g_B(x,\ddzi)}\ge\sqrt{\frac{n+1}b}(\frac ab)^{n}.$$  Since 
$a\leqslant \sqrt{g_K(x,\ddzi)}\leqslant b$ from Ahlfors Schwarz Lemma,
 we conclude that
$$\frac{2\pi}{a^3}(\frac {2b}a)^n\sqrt{g_K(x,\ddzi)}\ge \sqrt{g_B(x,\ddzi)}\ge\sqrt{\frac{n+1}{b^2}}(\frac ab)^{n}\sqrt{g_K(x,\ddzi)}.$$ 

As  $a^2\leqslant g_K(x,V),$ $g_K$ is non-degenerate on $M.$  The earlier 
argument estimating $g_B$ by $g_K$ from below then implies that $g_B$ is
non-degenerate.  Hence $g_B$ is a K\"ahler metric. 

We prove now that $g_B$ is complete.  If $g_B$ is incomplete, it follows
that there is a geodesic $\gamma$ of finite length $l$ from a fixed point $x_o\in M$
approaching to a point $y$ on $\pd M.$  In particular, given any preassigned
number $\epsilon>0,$ we can choose a point $z$ on $\gamma$ so that the distance 
$d_B(z,y)=\lim_{w\rightarrow y}d_B(z,w)\leqslant\epsilon.$  On the other hand,
the above discussions relating $g_B$
to $g_K$ actually shows that the
distance $d_B(z,\pd B_a(x))$ with respect to the Bergman metric is at least
$a\cdot k_1.$  This clearly leads to a contradiction by choosing 
$\epsilon<a\cdot k_1.$  Hence $g_B$ is complete.

\bs
\begin{lemm}
$M$ is a pseudoconvex domain. 
\end{lemm}

\ni{\bf Proof}
Fix a realization of $M$ as a bounded domain $\Omega$ in $\bC^n.$
From the previous lemma, the Bergman metric $\sqrt{-1}\pd\op\log K_{B,\Omega}$ is positive definite, here 
$K_{B,\Omega}=\sum_i|f_i(z)|^2$ is the potential of the Bergman metric
on $\Omega$ expressed in terms of  a unitary basis $\{f_i\}$ of 
the space of $L^2$-holomorphic
functions on $\Omega.$  
Clearly, $K_{B,\Omega}$ is a strictly plurisubharmonic function on $M.$
We need only to prove that $K_{B,\Omega}$ blows up
along any sequence of points approaching the boundary of $\Omega.$

For a point $x\in \Omega\cong M,$ let us still use the notation $\varphi_x$ for
the uniformizing coordinate charts for $x$ as defined in the Introduction.
Let $K_{B,\varphi_x(\Omega)}$ be the potential of the Bergman metric on 
$\varphi_x(\Omega).$  Clearly in terms of the Jacobian of the transition
functions,
$$K_{B,\Omega}=K_{B,\varphi_x(\Omega)}|J(\varphi_x)|^2.$$
From the proof of Lemma 1, we know that 
\begin{eqnarray*}
K_{B,\varphi_x(\Omega)}(y,y)&=&\sup_{f\in H^0_{(2)}(\varphi_x(\Omega)),\Vert f\Vert_{\varphi_x(\Omega))}=1}|f(y)|^2\\
&\ge&\frac 1{\int_{B_b(0)}1}
\end{eqnarray*}
is bounded from below.  Hence it suffices for us to prove that 
$|J(\varphi_x)|$ blows up for $x$ approaching $\pd \Omega.$

Recall from definition that $\varphi_x^{-1}:\varphi_x(\Omega)\rightarrow \Omega$ is a biholomorphism and
$\varphi_x^{-1}(0)=x.$  We claim that as $x\rightarrow\pd\Omega,$ 
the smallest eigenvalue $\mu_x$ of the Jacobian matrix $J(\varphi_x^{-1})\vert_0$ at $0$ approaches to $0.$ 

To prove the claim, we assume for the sake of proof by contradiction that there exists a sequence of points $x_i\in \Omega$ with
Euclidean distance $d(x_i,\pd\Omega)=\epsilon_i\rightarrow 0$ but $\mu_{x_i}\geqslant c_1$ for some constant $c_1>0.$
First of all, we observe by applying the generalized Cauchy estimates to $\varphi_x^{-1}$ on $B_a(0)$ that 
every derivative of $\varphi_x^{-1}$ is bounded from above by some constant independent of $x.$  In particular, all second derivatives of $\varphi_x^{-1}$
with respect to the coordinate vectors on $B_a(0)$ are bounded from above by a constant $c_2>0.$
Let now $\ell_i$ be a line segment in $\Omega$ realizing the Euclidean distance between $x_i$ and $\pd\Omega,$ so
that for $y_i\in \ell_i\cap \pd\Omega,$ $d(x_i,y_i)=\epsilon_i.$  The complexification of $\ell_i$ is
a complex line $\ell_{i,\bC}$ intersecting $\Omega.$  After a linear change of coordinate, we may assume that
$\ell_{i,\bC}$ is defined by $\zeta_2=\cdots=\zeta_n=0.$  We may also assume that $\zeta_1=0$ at $x_i.$  Writing $\zeta=t+i u$ in terms of real and imaginary part, we may assume without loss of generality
that $\ell_i$ lies on the
real axis defined by $u=0$ and parametrized by $t$ for $0\leqslant t\leqslant \epsilon_i.$   Hence the end point
on $\pd\Omega$ is given by $\zeta_i=0,$ $i\geqslant 2,$ and $\zeta_1(y_i)=\epsilon_i.$

As $\varphi_{x_i}$ is a biholomorphism, the image $\tell_i:=\varphi_{x_i}(\ell_i\cap \Omega)$ is a real curve on $\varphi_{x_i}(\Omega)$ with $0$ as an endpoint.  Assume that
$\tell_i$ is parametrized by a unit speed parameter $s$ so that $\tell_i(0)=0$ on $\varphi_{x_i}(\Omega).$   Since $\varphi_{x_i}^{-1}(B_a(0))$ intersects $\ell_i$ on $\Omega,$ we know that the length of $\tell_i$ in $\varphi_{x_i}$ is greater than $a.$
Let $r=\min(a,\frac{c_1}{4c_2}).$   Denote by $\lambda(s)$ the minimal eigenvalue of $J(\varphi_x^{-1})$ at $\tell(s).$
$\lambda(0)=\mu_{x_i}$ from definition.
From the Mean Value Theorem in calculus, it is clear that for $0\leqslant s\leqslant r,$
the minimal eigenvalue at $s$ satisfies
$$\lambda(s)\geqslant \lambda(0)-c_2s\geqslant c_1-c_2r\geqslant \frac {c_1}2.$$
It follows that $\frac{dt}{ds}\geqslant \frac {c_1}2$ for $0\leqslant s\leqslant r.$  Hence the length of $\ell_i$ is at
least $\frac {c_1}2 r,$ a constant independent of $x_i.$  Clearly this contradicts the assumption that the length of
$\ell_i,$ which is $d(x_i,\pd\Omega)$, is $\epsilon_i$ and $\epsilon_i\rightarrow 0$ as $i\rightarrow\infty.$  The claim is proved.

As mentioned above, each 
eigenvalue of the Jacobian of $\varphi_x^{-1}$ 
is bounded from above by a constant $c_3$ for all points $x\in \Omega.$
Moreover, the smallest of them approaches to $0$ as $x\rightarrow\pd\Omega.$
Since the determinant $|J(\varphi_x^{-1})|$ is just the product of all the eigenvalues of the $J(\varphi_x^{-1}),$
we conclude that
$|J(\varphi_x^{-1}|\rightarrow 0$ as $x\rightarrow\pd\Omega$.  Hence 
  $|J(\varphi_x)|$ tends to 
$\infty$ as $x\rightarrow \pd\Omega.$

\bs
\ni{\bf Remark}  It was pointed out by the referee that the argument essentially showed that the trace of the Bergman kernel 
$K(x,x)$ of a uniform squeezing domain was
bounded from below by $c/d,$ where $d=d(x,\pd\Omega)$ is the Euclidean distance to the boundary of $\Omega$ 
and $c$ is a constant.   Later on we will see that
as a consequence of Theorem 1 and 2, the estimates can be improved to $c/(d^2(-\log d)^2)$ as stated in Corollary 3.

\bs
We may now complete the proof of Theorem 1.

\ms
\ni{\bf Proof of Theorem 1}  (a) and (b) follow from Lemma 1 and Lemma 2.
(c) follows from the work of Cheng-Yau and Mok-Yau on K\"ahler-Einstein metrics
(cf. [MY]).

\bs
\begin{center}
{\bf \S3 Metric properties}
\end{center}

We say that we metrics $g_1$ and $g_2$ are equivalent or quasi-isometric on a domain $\Delta$,
denoted by $g_1\sim g_2,$ if there exists a constant $c>0$ such that $\frac1c g_2\leqslant g_1\leqslant cg_2.$

\begin{prop} The invariant metrics on a uniformly squeezing domain satisfy $g_C\sim g_K\sim g_B\sim g_{KE}.$
More precisely, 
\begin{eqnarray*}
\frac ab g_K\leqslant&g_C&\leqslant g_K,\\
\frac ab g_K\leqslant&g_B&\leqslant[\frac{2\pi}{a^3}(\frac {2b}a)^n]^2g_K,\\
\frac{a^2}{b^2n}g_K\leqslant&g_{KE}&\leqslant\frac {b^{4n-2}n^{n-1}}{a^{2n-2}}g_K.
\end{eqnarray*}
\end{prop}

\ni{\bf Proof} Since the proof is very similar to the proof of Theorem in [Y3],
we would just give a brief outline.

 It follows from Ahlfors Schwarz Lemma that
$g_C\leqslant g_K$.  On the other hand, from definition
of $g_K$ and $g_C$ and the inclusions $B^n_a(x)\subset \varphi(M)\subset B^n_b(x)$, we conclude for any tangent vector $v\in T_vM$ that $\sqrt{g_K(x,v)}\leqslant \frac1a$
and $\sqrt{g_C(x,v)}\ge b.$  Hence $g_B\leqslant \frac ba g_C.$  Hence 
$g_B\geqslant \frac ab$ from the above discussions.

The upper bound of $g_B$ by $g_K$ is already given in the proof of Lemma 1.
On the other hand, as observed by Look and Hahn (cf. [H]),
it follows by expressing $g_B$ in terms of extremal functions that
$g_K\ge g_C.$

To compare $g_{KE}$ and $g_K,$ we normalized
the Poincar\'e metric on $B_r(0)$ so that the potential is $\log(r^2-|z|^2).$ 
The resulting metric is a K\"ahler-Einstein metric of Ricci curvature $-\frac{2n+2}{r^2}$ with
constant holomorphic sectional curvature $-4<0$.  Then 
$g_{KE}^{B_r}(v,\ov)=\frac1{r^2}=g_{K}^{B_r}(x,v).$
It follows from definition that on $B_r,$ 
$g_K^{B_r}(0,v)=g_C^{B_r}(0,v)=[\frac1r]^2=g_{KE}^{B_r}(v,\ov).$
$B_{a}^{n}\subset M_S\subset B_{b}^{n}.$  Let us denote the volume form of
$g$ by $\mu(g).$
Applying Schwarz Lemma of Mok-Yau [MY]
to the first inclusion with respect to the K\"ahler-Einstein
metrics $g_{KE}^{B_{a}^{n}}$ and $g_{KE}^{M}$ on $B_{a}^{n}$
and $M$ of Ricci curvature $-\frac{2n+2}{(a)^2}$ and $-(2n+2)$ respectively, we get 
\begin{eqnarray*}
\mu(g_{KE}^{M})&\leqslant&a^{2n}\mu(g_{KE}^{B_{a}^{n}})=a^{2n}\mu(g_{K}^{B_{a}^{n}})\\
&=&b^{2n}\mu(g_{K}^{B_{b}^{n}})\leqslant b^{2n}\mu(g_{K}^{M})
\end{eqnarray*}
Applying Schwarz Lemma of [R] to $g_{KE}^{M}$ which has constant Ricci
curvature $-(2n+2)$ and 
$g_{KE}^{B_{b}^{n}}$ which has constant 
holomorphic sectional curvature 
$-4,$ we conclude that 
$$
g_{KE}^{M}\ge\frac1n g_{KE}^{B_{b}^n}=\frac1ng_{K}^{B_{b}^n}\ge 
\frac{a^2}{nb^2}g_{K}^{B_{a}^n}\ge
\frac{a^2}{nb^2}g_{K}^{M}.$$
Let $\nu_i>0, i=1,\dots, n$ be the eigenvalues of 
$g_{KE}^{M}$ with respect to $g_{K}^{M}.$  We conclude from the second estimate that 
$\nu_i\ge \frac{a^2}{nb^2}$ for all $i,$ and from the first statement that
$\prod_{i=1}^n\nu_i\leqslant b^{2n}.$  It follows that 
$\frac{b^{4n-2}n^{n-1}}{a^{2n-2}}\ge\nu_i\ge \frac{a^2}{b^2n}.$  Hence
$$(\frac{b}{a})^{2n-2}n^{n-1}g_K\ge g_{KE}\ge \frac {a^2}{b^2n}g_K.$$
This concludes the proof of the proposition.

\begin{prop}
(a). There exists a constant $c^{g_{KE}}_N$ depending only on the order of differentiation $N$
such that 
$\Vert \td^{g_{KE}}_{i_1,\cdots i_N}R_{g_{KE}}\Vert\leqslant c^{g_{KE}}_N$ for any covariant derivatives $\td^{g_{KE}}_{i_1,\cdots i_N}.$
Consequently, the curvature tensor of $g_{KE}$ and
any order of covariant derivatives of the curvature tensor is bounded by a
uniform constant.  Furthermore, the injectivity radius of $g_{KE}$ 
is bounded uniformly from
below on $M$.\\
(b). The same conclusion is true for the Bergman metric $g_B.$
\end{prop}

\ni{\bf Proof}
(a). We denote $g_{KE}$ simply by $g$ in this part of proof.
 For each fixed point $x\in M,$ there exists a uniformizing squeezing
coordinate system given by $B_a(0)\subset\varphi(M)\subset B_b(0)\subset\bC^n,$ where 
$\varphi(x)=0.$  We would derive our estimates on such coordinate neighborhoods.
By a unitary change of coordinates, we may assume that $g_{i\oj}(x)$ is diagonal 
at $x=0.$
Furthermore, from Lemma 3, we know that on $B_{\frac a2}(0),$ 
$g_{KE}\sim g_{K}\sim g_o.$
Hence in measuring the magnitude of a derivative
 with respect to $g=g_{KE}$, 
it is up to some uniform constant the same as
measuring with respect to the Euclidean metric $g_o$.  We need the following
technical estimates.

\begin{lemm}
Let $g=g_{KE}$ be the K\"ahler-Einstein metric on a domain $M$ with uniform
squeezing properties.  Then all the covariant derivatives of the coordinate vector fields
in terms of the uniform squeezing coordinate systems are uniformly bounded on $M$ by a constant 
depending on the order of differentiations.
\end{lemm}

\ni{\bf Proof}  In terms of the uniform squeezing coordinate system, Lemma 3 is equivalent to the
boundedness of any order of derivatives of the metric coefficient $g_{i\oj}$ with respect to the coordinate
vectors.

In the following, we denote by $c_i, c_{k,m}$ and $c'_{k,m}$ constants which are
independent of $x\in M.$
The K\"ahler metric satisfies Einstein equation
\begin{eqnarray}
\pd_i\op_j\log|\det(g)|&=&cg_{i\oj}
\end{eqnarray} on $B_a(x)$ for $x\in M.$

Note that $g_{i\oj}$ coming from solution of Monge-Amp\`ere equation 
is smooth from the standard results in K\"ahler-Einstein equations
(cf. [Au], chapter 7).  In fact, it would also follow from Proposition
2 together with Theorem 3.56 of [Au].
Taking trace with respect to the Euclidean metric, we get
\begin{eqnarray}
\Delta_o\log|\det(g)|=cg_o^{i\oj}g_{i\oj}
\end{eqnarray}
on $B_a(x).$

Recall also from Proposition 3 
that on $B_{a/2}(x),$ 
$g_{i\oj}$ is quasi-isometric to the Euclidean metric $|(g_o)_{i\oj}|=\delta_{ij},$
where $\delta_{ij}$ are Kronecker's delta.

In terms of the usual notion used in [GT], let us denote by $H_k=W_{k,2}$ and 
$C_{k,\alpha}$ the spaces of functions on $B_a$ which are
bounded with respect to the Sobolev norm and H\"older norm on $B_a(x)$ 
respectively.  Applying Calderon-Zymund's estimates in Proposition 2 to the Einstein equation,
we get 
$$\Vert\log|\det(g)|\Vert_{H_2}\leqslant c_1[\Vert\log|\det(g)|\Vert_{H_0}
+\Vert g_o^{i\oj}g_{i\oj}\Vert_{H_0}]\leqslant c_2,$$
here again we used Proposition 3.

Observe that on each point $y\in B_a(x),$ there exists a unitary matrix $A_y$
such that $A_yg\oA_y^{t}$ is a diagonal matrix.  As $|\det A_y|=1$ from 
definition, we may assume that $g$ is diagonal at $y$ for our computation
involving $\det(g).$  Hence $|\det(g)|=\prod_{i=1}^ng_{i\oi}$
at $y$.
The Einstein equation gives rise
to
$\pd_i\op_i\log|\det(g)|=cg_{i\oi}.$  Let $X$ be any vector field of unit
length coming from
linear combination of coordinate vector fields.  Let
$D_X$
denote derivative in the direction of $X$.  Then by applying $D_X$ to the
above equation
$$\pd_i\op_iD_X\log|\det(g)|=cD_Xg_{i\oi}=cg_{i\oi}D_X\log|(g_{i\oi})|.$$
Taking the trace by $g$ and summing over all $i=1,\dots,n,$ we get
$$\sum_{i=1}^ng^{i\oi}(y)\pd_i\op_iD_X\log|\det(g)|=c\sum_{i=1}^n D_X\log|(g_{i\oi})(y)|
=D_X\log|\det(g)|(y).$$
We obtain $\Delta_gD_X\log|\det(g)|=D_X\log|\det(g)|(y).$
As $g_{i\oj}$ is uniformly quasi-isometric to $(g_o)_{i\oj},$
 the same Schauder estimate allow us to conclude that
$$\Vert D_X\log|\det(g)|\Vert_{H_2}\leqslant c_3[\Vert D_X\log|\det(g)|\Vert_{H_0}]
\leqslant c_4$$ after applying the earlier bound on $\log|\det(g)|\Vert_{H_2}.$
As $X$ is arbitrary, this implies 
$$\Vert \log|\det(g)|\Vert_{H_3}
\leqslant c_4.$$  Clearly the bootstrapping argument implies that for each positive
integer $m$, there exists a constant $c_m$ independent of $x$ such that
$$\Vert \log|\det(g)|\Vert_{H_m}
\leqslant c_m.$$

Applying the Sobolev Estimates to $B_{a/2}(x),$ we conclude that 
$$\Vert \log|\det(g)|\Vert_{C_{k,m}}\leqslant c_{k,m}$$
for some constant $c_{k,m}$ independent of $x.$  Applying to equation (0.2),
we conclude that
 $$\Vert g_{i\oj}\Vert_{C_{k,m}}\leqslant c'_{k,m}$$ for some constant
$c'_{k,m}$ independent of $x.$  This concludes the proof of the Lemma.

\bs
We may now complete the proof of (a).  We compute in terms of the coordinate vectors in the uniform squeezing coordinate system.
A covariant derivative of the curvature tensor is a 
Euclidean derivative modified by an addition term coming from Christopher
symbol, which are linear combinations of the derivatives of the metric tensor with respect to coordinate
vectors.
It follows easily from Lemma 3 and by induction that the norm of an $n$-th order derivative
of the metric tensor with respect to the coordinate vectors is bounded by a constant depending on $n$ but 
 independent of $x.$  Hence the first 
statement of
(a) follows.

Since the curvature tensor of $g_{KE}$ appears as sum of some second order
and first order derivatives of the metric tensor with respect to coordinate vectors, 
it is clear that the curvature
tensor is bounded uniformly.  Similarly, any order $N$ covariant derivatives
of the curvature are linear combination of expressions involving up to $N+2$
order of derivatives of the metric tensor with respect to the coordinate vectors.  We conclude that any such derivative
is bounded by a constant depending on $N$. 

From the boundedness in curvature, we conclude immediately that 
the conjugate radius
is bounded from below by an absolute constant.  Furthermore we note that
 there exists a $\epsilon>0$
such that a
geodesic loop $l$ of length less than $\epsilon$ based at $x$ does not
exist.   Suppose on the contrary such a geodesic exists and the incoming
and outcoming geodesic segment span an angle $\theta$, which has to be
positive, at $x.$
Then we may find two points $P_1, P_2$ on $l$ near cut-locus of $x$
such that 
the distance of their preimages on the
tangent space with induced metric by the exponential map at $x$ is bounded
from below by $\theta\epsilon$, but clearly not on $M$.  
This is clearly a contradiction for $\epsilon$ sufficiently small and
the fact that the metric is quasi-isometric to the Euclidean one.

\ms
\ni(b). To consider the derivatives of the Bergman metric, again we consider the
uniformizing squeezing coordinates and let  
$K_x(z,w)=\sum_if_{i}(z)\overline{f_{i}(w)}$ be the coefficient of the
Bergman kernel on $\varphi_x(M)$
which is holomorphic in z and conjugate holomorphic in $w.$ Let $\oM$ be the set $M$ equipped with the conjugate 
complex structure.
Writing $w_i=\ou_i,$ we conclude that $K_x(z,\bar u)$
is holomorphic on $M\times \oM$ with respect to the 
complex structures on $M$ and $\oM$ respectively.
The restriction $K(z,w)$ to
$w=z$ is precisely the potential for the Bergman metric.

Let $A=B_{\frac a2}\times B_{\frac a2}.$
Let $D$ be a differential operator involving compositions of the
coordinate derivatives.
By Generalized Cauchy Inequality, it follows easily that
all the higher derivatives 
$|[D\ddzi\ddbwj \log{K}](z_o,\ou_o)|$ of the metric 
at the origin are controlled up to a constant
depending on $D$ by $| {K}(z,\ou)|$ for $(z,u)$ lying on 
 the boundary $\pd(A).$  Clearly $|K_x(z,\ou)|^2\leqslant K_x(z,z)K_x(\ou,\ou)$
by the Cauchy-Schwarz Inequality for $(z,u)\in\pd A.$  

In terms of the
peak function $f_z$ at $z\in \varphi(M)$ mentioned before, we obtain
\begin{eqnarray*}
K_x(z,z)&=&|f_z|^2\leqslant \frac1{\vol(B_{\frac a2}(z))}\int_{\vol(B_{\frac a2}(z))}|f_z|^2\\
&\leqslant&\frac1{\vol(B_{\frac a2})},
\end{eqnarray*}
since the $L^2$-norm of $f_z$ is $1.$  The same bound is applicable to
$K_x(\ou,\ou).$  Restricting to the twisted diagonal given
by $u=\bar z,$ it follows immediately that the curvature
tensor and all their derivatives are bounded
with respect to the Euclidean metric on $B_{\frac a2}(z)\subset M$.
As $g_B$ is uniformly quasi-isometric to $g_o$ on $\vol(B_{\frac a2}),$
  we conclude that all the derivatives
of the curvature are uniformly bounded for the Bergman metric.  As in part (a),
the finishes the proof of (b).

\bs

\begin{prop}
(a). $M$ is K\"ahler-hyperbolic with respect to $g_{KE}.$\\
(b).  The same is true for $g_B.$
\end{prop}

The proof of the proposition depends on the following lemma.

\begin{lemm}
Let $g=g_{KE}.$
Fix $x,y\in M$.  Let $W$ be a $(1,0)$-vector at $0\in \varphi_y(M)$
Let $
\varphi_{y,x}:\varphi_y(M)\rightarrow \varphi_x(M)$ be the
biholomorphic mapping given by $\varphi_x\circ\varphi_y^{-1}.$  
Then
$\frac{|\pd_W\log|J(\varphi_{y,x})|^{2}|}{\sqrt{g(W,\oW)}}(0)\leqslant C$
for some constant $C$ independent of $x$ and $y.$  
\end{lemm}

\ni{\bf Proof}  Let $z^i, w^j, i,j=1,\dots,n$ be the local coordinates
on $\varphi_y(M)$ and $\varphi_x(M)$ respectively.
For simplicity, we also denote $\varphi_{y,x}$ by
$\varphi$ since $x$ and $y$ are fixed in this proof.  Let us choose the
coordinate at $\varphi_x(M)$ such that $g(\ddwi,\oddwj)$ is diagonal at
$\varphi(0).$

It suffices for us to show that for each $k=1,\dots,n,$
$\frac{|\frac{\pd}{\pd z^k}\log|J(\varphi_{y,x})|^{2}|}{\sqrt{g(\ddzk,\oddzk)}}(0)\leqslant C.$  Clearly,
$$|J(\varphi)(z)|^2=\frac{\det(g(\ddzi,\oddzj))(z)}{\det(g(\ddwi,\oddwj))(\varphi(z))}.$$
Hence
\begin{eqnarray*}
&&\frac{|\frac{\pd}{\pd z^k}\log|J(\varphi_{y,x})|^{2}|}{\sqrt{g(\ddzk,\oddzk)}}(0)
\leqslant
[\frac1{\sqrt{g(\ddzk,\oddzk)}}\frac{\pd}{\pd z^k}|\log|\det(g(\ddzi,\oddzj))||]_{z=0}\\
&&+[\frac1{\sqrt{g(\ddzk,\oddzk)}}\frac{\pd}{\pd z^k}|\log|\det(g(\ddwi,\oddwj))(\varphi(z))||]_{z=0}.
\end{eqnarray*}

For the first term,  applying Schauder type estimates to the chart $\varphi_y(M)$
and using the fact that $g=g_{KE}$ 
is K\"ahler-Einstein, we conclude  as in the proof of Proposition 3 that
$$[|\frac1{\sqrt{g(\ddzk,\oddzk)}}\frac{\pd}{\pd z^k}|\log|\det(g(\ddzi,\oddzj))(z)||]
_{z=0}<c_1,$$
Note that quasi-isometry of $g=g_{KE}$ and the Euclidean metric is used here.

For the second term, we rewrite $w=\varphi(z)$ and use Chain rule to
rewrite

\begin{eqnarray*}
&&[\frac1{\sqrt{g(\ddzk,\oddzk)}}\frac{\pd}{\pd z^k}|\log|\det(g(\ddwi,\oddwj))(\varphi(z))||]_{z=0}
\\
&=&[\frac1{\sqrt{\sum_{i,j}g(\ddwi,\oddwj)\frac{\pd w^i}{\pd z^k}
\overline{\frac{\pd w^j}{\pd z^k}}}
}
|\sum_l\frac{\pd w^l}{\pd z^k}\frac{\pd}{\pd w^l}|\log|\det(g(\ddwi,\oddwj))(\varphi(z))|]_{z=0}\\
&=&\frac1{\sqrt{\sum_{i}g(\ddwi,\oddwi)|\frac{\pd w^i}{\pd z^k}|^2(0)}}
[|\sum_l\frac{\pd w^l}{\pd z^k}(0)\frac{\pd}{\pd w^l}|\log|\det(g(\ddwi,\oddwj))(w)|]_{w=\varphi(0)}.
\end{eqnarray*}

Applying now Schauder's estimate to $\varphi_x(M)$ and using the fact that
$g$ is K\"ahler-Einstein and $g$ is quasi-isometric to the Euclidean metric
on $\varphi_x(M)$ again, we conclude that $|\frac{\pd}{\pd w^l}|\log|\det(g(\ddwi,\oddwj))(w)|_{w=\varphi(0)}\leqslant c_2.$  Hence 
$$\frac1
{\sqrt{\sum_ig(\ddwi,\oddwi)(\varphi(0))|\frac{\pd w^i}{\pd z^k}|^2(0)}}
|\frac{\pd w^l}{\pd z^k}(0)
\frac{\pd}{\pd w^l}|\log|\det(g(\ddwi,\oddwj))(\varphi(0))|
\leqslant c_3$$
 and we conclude that 
$$\frac{|\frac{\pd}{\pd z^k}\log|J(\varphi_{y,x})|^{2}|}{\sqrt{g(\ddzk,\oddzk)}}(0)
\leqslant c_1+c_3.$$
As $k$ is arbitrary, this concludes the proof of the Lemma.

\bs
\ni{\bf Remark} Lemma 4 gives a proof of the part of Proposition 5 of [Y4] about boundedness of the first derivative
of the Jacobian, showing that it corresponds
to general properties of uniform squeezing domains.  The last sentence in Proposition 5 of [Y4] that $|J(\Phi_{y,x})|=1$
is incorrect and was not used in the subsequent arguments in [Y4].

\bs
\ni{\bf Proof of Proposition 5}
(a) Fix a point $x\in M$ and consider the uniform squeezing coordinates
$\varphi_x$
associated to it, so that $B_a(0)\subset\varphi_x(M)\subset B_b(0).$
The K\"ahler form $\omega_{KE}$ associated to the K\"ahler-Einstein metric
$g_{KE,x}$ satisfies
$$\sqrt{-1}\pd\op\log\det g_{KE}=c\omega_{KE}$$
for some negative constant $c.$  Left hand side of the above expression is
independent of the particular coordinate $\varphi_x$ that we are using.
Note that the determinant $\det(g_{KE}(\ddzi,\oddzj))$ depends on
the local coordinates $\varphi_x(M).$  We denote the quantity by $\det(g_{KE,x}).$
In this way we may regard
$h_x=\sqrt{-1}\op\log\det g_{KE,x}$ as the potential one form to satisfy
$dh_x=\omega,$ here note that $h_x$ depends on our fixed base point $x.$
 Again, on $B_{\frac a2}(0),$ the equation is
$$\Delta_{g_{KE,x}}\log\det g_{KE,x}=c.$$
Since $g_{KE,x}\sim g_K\sim g_o$ on $B_{\frac a2}(0),$ the above equation is
a strongly elliptic equation with uniformly bounded coefficients.  It follows
from Proposition 3 that we have a bound
$|d\log\det g_{KE,x}|\leqslant C$ for some uniform constant $C.$
We conclude that $|\sqrt{-1}\op\log\det g_{KE,x}|<C$ and hence
$$|h_x|_{g_{KE,x}}=|\sqrt{-1}\op\log\det g_{KE,x}(y)|_{g_{KE,x}}<C_1$$ for some uniform constant 
$C_1$ and all $y\in B_{\frac a2}(0).$

Now we need to worry about points $y\in \varphi_x(M)-  
B_{\frac a2}(0).$  For such cases, we consider the uniform squeezing 
coordinate $\varphi_{y}$ as well, here we identify $y$ with $\varphi_x^{-1}(y)$
to simplify our notations.
$\varphi_{y,x}=\varphi_y\circ \varphi_x^{-1}$ is the biholomorphism from $\varphi_x(M)$ to $\varphi_y(M).$
We have correspondingly $\det(g_{KE,x}(z))=\det(g_{KE,y}(w))|J(\varphi_{y,x})|^2.$ 
Let $Y\in T_0(\varphi_y(M))$ and 
$X=(\varphi_{y,x})_*Y\in T_{\varphi_{y,x}(0)}(\varphi_x(M).$
Then
\begin{eqnarray*}
h_x(X)&=&\sqrt{-1}\pd_X\log \det(g_{KE,x})=\sqrt{-1}\pd_Y\log \det(g_{KE,y})+\sqrt{-1}\pd_X\log|J(\varphi_{y,x})|^2\\
&=&h_y+\sqrt{-1}\pd_X\log|J(\varphi_{y,x})|^2.
\end{eqnarray*}

Clearly
it follows from the earlier argument that $|h_y|_{g_{KE,y}}<c_2$
for some uniform constant $c_2.$  Lemma 4 also shows that
$|\pd_X\log|J(\varphi_{y,x})||^2_{g_{KE,y}}<c_3$.  Proposition 3 for $g_{KE}$
now follows from combining these two estimates.

\ms
\ni(b)  
Fix our point $x\in M$ as before.  The K\"ahler form of $g_B$ can
be written as $\omega_B=d\eta_x,$ where
$\eta_x=-\sqrt{-1}K_{x}^{-1}\pd K_{x}$ and $K_x\sum_{j}|f_j(z)|^2$ by taking over unitary
basis of holomorphic functions on the Hilbert space of 
$L^2$-holomorphic functions on $\varphi_x(M).$  Proposition 3b above implies that
$|\eta_x(y)|_{g_B}<c_4$ for some constant $c_4$ and every point $y\in B_{\frac a2}(0)\subset \varphi_y(M).$  

For a point $y\in \varphi_x(M)-  
B_{\frac a2}(0),$ we note that the potentials of $g_B$ satisfies
$K_x=K_y|J(\varphi_{y,x})|^2.$  Hence
$$\eta_x(X)=
\eta_y(Y)+\sqrt{-1}\pd_X\log|J(\varphi_{y,x})|^2$$
similar to the derivation in (a).
From the previous paragraph, $|\eta_y|_{g_B}<c_5$ for some uniform
constant $c_5>0.$  From Proposition 3b again,
$|\pd_X\log|J(\varphi_{y,x})(y)|^2_{g_B}\leqslant c_6|\pd_X\log|J(\varphi_{y,x})(y)|^2_{g_o}|\leqslant c_7$ for some constant $c_7>0.$  (b) follows
 by combining the previous two estimates.

\begin{lemm}
Let $x$ be a fixed point on $M$.
Expressed in terms of the uniform squeezing coordinate system $\varphi_x(M),$
$|\det(g_{KE})|^{-\alpha}$ for $\alpha$ sufficiently small
is a bounded plurisubharmonic exhaustion function
on $M.$
\end{lemm}

\ni{\bf Proof}
Denote by $|g_{KE}|=\det(g_{KE})$ the determinant of $g_{KE}$ in local coordinates.
Direct computation yields
\begin{eqnarray*}
&&\sqrt{-1}\pd\op(-|g_{KE}|^{-\alpha})\\
&=&-\alpha(\alpha+1)|g_{KE}|^{-\alpha-2}\sqrt{-1}\pd |g_{KE}|\wedge \op |g_{KE}|+\alpha |g_{KE}|^{-\alpha-1}\sqrt{-1}\pd\op |g_{KE}|\\
&=&\alpha|g_{KE}|^{-\alpha}\sqrt{-1}[(\alpha+1)\pd\op(\log |g_{KE}|)-\pd\log |g_{KE}|\wedge\op\log |g_{KE}|].
\end{eqnarray*}
Applying Proposition 5
and noting that $\log |g_{KE}|$ is up to a constant the potential of $g_{KE}$
on any realization of $M$ as a bounded domain, $|\sqrt{-1}\pd\log |g_{KE}|\wedge\op\log |g_{KE})||\leqslant
c\sqrt{-1}\pd\op\log(|g_{KE}|)$ for some constant $c>0.$  It suffices for
us to choose $\alpha>\frac1c-1$ to conclude the proof of the lemma.

\bs
\ni{\bf Proof of Theorem 2}
 (a) follows from Proposition 3.  (b), (c) and (d) follows from Proposition 4.   (e) follows from Proposition 5. (f) follows from Lemma 5.

\bs
\ni{\bf Remark} 
We remark that finiteness in volume of $g_{KE}$ is in fact equivalent to
the quasiprojectiveness of our $M$.  One direction is proved in the above
corollary.  For the other direction , as assume that
$M_1=M/\Gamma$ is a quasi-projective manifold.
We may assume that $M_1=\oM_1-D$ for some normal crossing divisor $D$
after resolution of singularities if necessary.
Hence neighborhoods of $D$ in $M_1$ are covered by union of open
sets $U_i$ of the form $\Delta_1^a\times (\Delta_1^*)^b,$
where $\Delta_1$ is a Poincar\'e disk of radius $1$ and $\Delta_1^*$ is a punctured
disk of radius $1$.

Equip each such $U_i=\Delta_1^a\times (\Delta_1^*)^b,$ we consider a smaller
open set $U_r=\Delta_r^a\times (\Delta_r^*)^b$ for $0<r<1$
with the restriction of the Poincar\'e metric $g_P$ on $U_i$ and apply the Schwarz 
Lemma
of Mok-Yau [MY] to the embeddings of the inclusion of $(U_i,g_P)$ into 
$(M_1,g_{KE}),$ we conclude easily that the volume of $(M_1,g_{KE})$ is
finite since the volume of each $(U_i,g_P)$ is finite.

\bs
\begin{center}
{\bf \S4 Examples}
\end{center}

We are going to show that examples in Proposition 1 do satisfy the uniform squeezing property.

\bs
\ni{\bf Proof of Proposition 1}
We will prove (a), (b), (d) and (e) first and leave the proof of (c) to the
end.\\
{\it (a) Bounded homogeneous domain $M$ in $\bC^n.$}

  Choose any point on
$M$ and translate the origin of $\bC^n$ to that point.  As it is bounded, it
is contained in a ball $B_b(x).$  Since $0$ lies in the interior of $M$,
there exists a ball of positive radius $a$ such that $B_a(0)\subset M.$
Hence $0\in B_a(0)\subset M \subset B_b(x)\subset \bC^n.$  Let $x$ be an
arbitrary point on $M$.  As $M$ is 
homogeneous, there exists a biholomorphism of $M$ moving $x$ to $0,$ 
here we realize $M$ as a fixed domain in $\bC^n.$  Hence balls of the
same radii provide uniform squeezing coordinates for
all $x\in M.$

\ms
\ni{\it (c) Bounded domains which cover a compact K\"ahler manifold}

Assume that $M\subset B_c(x_o)\subset \bC^n$ is a bounded domain covering a compact K\"ahler manifold $N$, where $x_o$ is a fixed point on $M$.
Let $A$ be a fundamental domain of $N$ in $M$.  For each point $x\in A,$ there 
exists a ball of radius $r_x$ such that $B_{r_x}(x)\subset M.$
Since $A$ is relatively compact, $r=\inf_{x\in A}r_x>0.$
It is then clear that $B_a(x)\subset M\subset B_b(x)$ for each $x\in A.$  Since each point
$y\in M$ can be mapped biholomorphically to some point in $A$ by the deck transformation
group, it is clear that we get a $(a,b)$ uniform squeezing coordinate.

\bs
\ni{\it (d) Teichm\"uller spaces $\cT_{g,n}$ of compact Riemann surfaces
of genus $g$ with $n$ punctures}

This is a consequence of Bers Embedding Theorem described as follows (cf. [Ga]).
Let $S$ be a Riemann surface of genus $g$ with $n$ punctures 
representing a point $x\in\cT_{g,n}.$  Denote by $\cT_S$ the Teichm\"uller
space based at $x.$
There exists an embedding $\Phi:\cT_S\rightarrow\bC^{N},$
so that $B_{\frac12}^{N}\subset \cT_S\subset B_{\frac32}^{N},$ 
where $\bC^{N}$ is identified with the space of holomorphic quadratic differentials
based at $S$ equipped with $L^\infty$ norm, and $\Phi(x)=0,$  where $N=3g-3+n.$

Hence the charts associated to Bers embedding provide us the uniform squeezing coordinates. 

\ms
\ni{\it (b) Bounded smooth strongly convex domains}

We give a step by step construction of the uniform squeezing coordinate systems. 

\ms
\ni(i) We observe the following fact.  Suppose $C_1=\pd B_a^1(x)$ and $C_2=\pd B_b^1(y)$ are two circles in $\bC$ of radii $a$ and $b$ meeting
tangentially at one point.  Assume that $B_a^1(x)\subset B_b^1(y).$  Let $w\in B_a^1(x)$ lying on the real line
joining $x$ and $y.$  Then there exists
a M\"obius $f$ mapping $C_1$ to itself, so that $f$ is holomorphic on $B_b^1(y),$ $f(w)=0$ and $f(C_2)\subset B_{2b}^1(0).$

To see this, we may assume that $x=0$ and $a=0$ by rescaling.  By a linear change of coordinates, we
may also assume that $y=-b+1$ lies on the real axis of $\bC.$  The fact follows by inspecting the explicit 
M\"obius transformation $z\rightarrow (z-w)/(1-z\ow).$

\ms
\ni(ii) We claim the following fact.  Suppose $B_a(x)\subset B_b(y)$ are two balls in $\bC^n$ and $\pd B_a(x)$
is tangential to $\pd B_b(y)$ at a point $q.$  Let $w\in B_a(x)$ lying on the real line joining $x$ and $y.$
Then there exists a M\"obius transformation $\psi$ of $B_a(x),$ 
so that $\psi$ is biholomorphic on $B_a(x),$  $\psi$ is holomorphic on $B_b(y),$ $\psi(w)=0$ and $\psi(B_b(y))\subset B_{2b}(0).$

To see this, after a linear change of coordinates, we may assume that the real line joining $x$ and $y$
is defined by $z_2=\cdots z_n=0$ and $Im(z_1)=0.$  As in (i), we may assume that $x=0$ by an affine change 
change of coordinate, and $a=1$ after rescaling.  Consider now the M\"obius transformation given by
$$\psi(z_1,\cdots,z_n)=(\frac{z_1-w}{1-z_1\ow},\frac{\sqrt{1-|w|^2}}{1-z_1\ow}z_2,\cdots,\frac{\sqrt{1-|w|^2}}{1- z_1\ow}z_n).$$
The same computation as in (i) establishes the claim.

\ms
\ni(iii) We now proceed to construct the uniform squeezing coordinate system.  We are
considering a $C^2$-strongly convex domain $M$ in $\bC^n.$
Let 
$p\in\pd M.$
Let $U_p$ be a neighborhood of $p$. For a point 
$q\in U'_p:=\pd M_p\cap U_p,$
let $N_{U'_p}(q)$
be the real line which is normal to $\pd M$ at $q$ with respect to the Euclidean metric.  As $\pd M$ is
$C^2$-smooth and $M$ is convex, there exist point $x_{p,q}, y_{p,q}\in N_{U'_p}(q),$ 
 and positive numbers $a_{p,q}$ and $b_{p,q}$ such that
both $\pd B_{a_{p,q}}(x_{p,q})$ and $\pd B_{b_{p,q}}(y_{p,q})$ 
are tangential to $\pd M'_p$ at $q$ and $B_{a_{p,q}}(x_{p,q})\subset M\subset B_{b_{p,q}}(y_{p,q}).$  

Replacing $U_p$ by a slightly smaller relatively
compact subset of itself if necessary, we may assume that 
$a_p=\liminf_{q\in U'_p}a_{p,q}>0$ and $b_p=2\limsup_{q\in U'_p}b_{p,q}<\infty.$  
Let $x'_{p,q}$ be the unique point on the normal line $N_{U'_p}(q)\cap M$ at a distance $a_p$ from $q.$
Let $V_p=\cup_{q\in U'_p} B_{a_p}(x'_{p,q}).$ 
From the above construction and from the claim in (ii), an $(a_p,b_p)$-uniform squeezing coordinate charts exists for
$V_p.$

The union $\cup_{p\in \pd M}V_p$ covers a neighborhood of $\pd M.$
From compactness of $\pd M,$ we can choose a finite number of points $p_1,\cdots,p_N$ on $\pd M$ such that $\cap_{i=1}^N V_{p_i}$ covers a neighborhood of $\pd M.$  Let $V_o$ be a relatively compact open subset of $M$ containing 
$M-\cap_{i=1}^N V_{p_i}$ so that $\{V_i\}_{i=0,\dots,N}$ gives a holomorphic
covering of $M.$  It is clear that there
exists 
$0<a_0<b_0$ such that for each point $z\in V_0,$ there exists a holomorphic
coordinate charts with $B_{a_0}(z)\subset M\subset B_{b_0}(z)$.
Let $a=\min(a_0,a_i, 1\leqslant i\leqslant N)$ and $b=\max(b_0,b_i, 1\leqslant i\leqslant N).$
It follows from our construction that the balls of radii $a$ and $b$ involved
form an $(a,b)$-uniform
squeezing coordinate system for $M.$  Hence strictly convex domain with $C^2$ 
boundary satisfies the uniform squeezing property.

 This concludes the proof of Proposition 1.

\bs

\bs
\begin{center}
{\bf \S5 Geometric consequences}
\end{center}

In this section, we give a proof of Corollary 1, Corollary 2 and Theorem 3 
 as applications
of Theorem 1 and 2.  A proof for Theorem 4 is also explained.

\bs
\ni{\bf Proof of Corollary 1}
(a) and (c) follows from the argument in [Gr].  (b) is already
 proved in [M], 
once we know that the manifold involved is K\"ahler-hyperbolic.  

\bs
\ni{\bf Proof of Corollary 2}
(a). It is given that $N=M/\Gamma$ is a compact complex manifold.  Since
we have already proved that there exists a K\"ahler-Einstein of negative
scaler curvature on $M$ which
descends to $N$,  it follows immediately that the canonical line 
bundle is ample and hence the variety is of general type.

\ms
\ni
(b). 
 From assumption, $M/\Gamma$ has finite volume with respect to the K\"ahler-Einstein
metric $g_{KE}.$
From Theorem 2a, $g_{KE}$ is complete K\"ahler,
with constant negative Ricci curvature and bounded Riemannian sectional
curvature.  Hence we may apply the results of [Y1], which relies on
the earlier results of Mok-Zhong [MZ], to conclude that $M_1$ is 
quasi-projective.

We note that  $h^i_{(2)}(M_1, lK_{M_1})=0$
for $l\ge2$
from $L^2$-estimates as in the corresponding proof of Kodaira's Vanishing Theorem.
Hence $h^i_{(2)}(M_1, lK_{M_1})=\chi(M, lK_{M}).$  It follows from 
a generalized form of the Riemann-Roch estimates of Demailly [De]
as proved in [NT] that
$$\chi(M_1,lK_{M_1})=\frac{l^n}{n!}\int_M c_1(K_M)^n+o(l^{N}).$$
Here we note that $c_1(K_{M})$ is positive definite from construction.
$L^2$-holomorphic sections of $lK_{M}$ extends as a holomorphic
section in $H^0_{\oM_1}(\oM_1,l(K_{\oM_1}+D)).$  This
follows from the fact that they
extend as $L^2$-sections and hence cannot have poles of order greater than
$1$ along any component of the compactifying divisor.  We conclude
that $\dim(H^0_{\oM_1}(\oM_1,l(K_{\oM_1}+D)))\ge cl^N$
and hence that $(\oM_1,D)$ is of log-general type.
This concludes the proof of Corollary 2.

\bs

\ni{\bf Proof of Theorem 3}
It follows from Proposition 1 that a bounded domain which
is the universal covering of a complex manifold is equipped with
a uniform squeezing coordinate.  Hence from Theorem 1, 
it is pseudoconvex and hence a Stein manifold.
Furthermore, it supports a complete K\"ahler-Einstein $g_{KE}$
metric
of negative scaler curvature.  Since $g_{KE}$ is invariant under 
biholomorphism and hence the deck transformations, it descends to
$N$.  Hence the canonical line bundle of $N$ is ample and the manifold
is of general type. 
We denote by $h^{p,q}_{(2),v}(M)$ the von Neumann dimension of
the space of $L^2$-harmonic $(p,q)$-forms
on $M$ with respect to the K\"ahler-Einstein metric (cf. [At]).
Corollary 1 implies that the von-Neumann dimension $h^{p,q}_{(2),v}(M)=0$ for
$p+q<n$ and $h^{n,0}_{(2,v)}(M)>0,$ which implies that the corresponding
Euler-Poincar\'e characteristics $(-1)^n\chi_{L^2,v}(M)>0.$  
From Atiyah's Covering Index Theorem,
$\chi_{L^2,v}(M)=\chi(M/\Gamma).$  Hence $(-1)^n\chi(M/\Gamma)>0.$

We may apply the same argument to the holomorphic line
bundle $2K$ on $M.$  We use $K$ to denote by the canonical line bundle on 
$N$ and $M.$  First of all $h^0_{L^2,v}(M,2K)>0$ by the usual 
$L^2$-estimates as used in [Y2], noting that $g_{KE}$ has
strictly negative Ricci curvature.  The same $L^2$-estimates implies that
$h^i_{L^2,v}(M,2K)=0$ for $i>0.$  Atiyah's Covering Index Theorem
implies that $\chi(N,2K)>0.$  On the other hand, from Kodaira's
Vanishing Theorem or $L^2$-estimates, we conclude that
$h^i(N,2K)=0$ for $i>0.$  Hence $h^0(N,2K)=\chi(N,2K)>0.$
This concludes the proof of Theorem 3. 

\bs
\ni{\bf Proof of Corollary 3}
From Theorem 1, we know that $\Omega$ is pseudoconvex.
It was proved by Mok-Yau in [MY] that a complete K\"ahler-Einstein metric of negative sectional
curvature exists and its volume form is bounded from below by $\frac1{d^2(-\log d)^2}$ with
respect to the Euclidean coordinates.  Corollary 3 follows from the proof of Theorem 2 as the Bergman
kernel is shown to be equivalent to the K\"ahler-Einstein volume form.  Note that both of them
transforms under a coordinate change by the same Jacobian determinant  as in the proof of Lemma 2.

\bs
\ni{\bf Proof of Corollary 4}
Stehl\'e has proved in [St] the result that a locally trivial holomorphic 
fiber space with hyperconvex fibers and Stein base is Stein.  Corollary 3
follows immediately from Theorem 2f.

\ms
\ni
{\bf Remark} There are many positive results to Serre's problem, including
the result of Siu [Si] when the fibers have trivial first Betti number,
the result of Mok [Mo] when the fibers are Riemann surfaces, and
the results of Diederich and Fornaess [DF].  In general
the problem has negative solution due to counterexamples such as the one given
by Skoda in [Sk].

\bs
\ni{\bf Proof of Theorem 4}
As explained in the last section, Bers Embedding gives rise to a uniform
squeezing coordinate system.   All the results of Theorem 4a-e follow from the earlier
results of this paper under the sole assumption that a uniform squeezing
coordinate system exists for $\cT_{g,n}$, which is provided by Bers Embedding.
Theorem 4f also follows from
Corollary 2 if we accept that $\cM_{g,n}$ is quasi-projective, which is
known classically by the well-known work of Baily for $n=0$ and Knudsen-Mumford
for the case of $n\neq0$ (cf. [KM]).

The exact formula for the Euler characteristic of $M$ has
already been obtained by Harer-Zagier [HZ].
We just remark that using K\"ahler hyperbolicity and Atiyah's Covering Theorem 
as in Theorem 3, we may prove that Euler-Poincar\'e characteristic of $M_1$ satisfies
$(-1)^n\chi(M_1)>0$ as well.  The only minor difference is that $M_1$ is
now non-compact.  However, it follows from Theorem 2 
that $(M_1, g_{KE})$ has finite geometry and hence the chopping
argument of Cheeger-Gromov [CG] shows that one may exhaust $M_1$ by
appropriate relatively compact sets so that the contribution from the 
boundary tends to $0$ as one takes the limit on the exhaustion.  It follows
that $(-1)^n\chi(M_1)>0.$

This concludes the proof of Theorem 4.  

\bs
\ni{\bf Remark} It is known that $g_K, g_C, g_B, g_{KE}, g_T$ and $g_M$ are quasi-isometric on $\cM_{g,n}$, where $g_T$ is the Teichm\"uller
metric and $g_M$ is a K\"ahler metric constructed in [Mc], as is shown in
[Y3].
It is also proved in [Mc] that any order of derivatives of $g_M$ 
is bounded as well.  Combining with Theorem 2c and the proof there, 
the difference $\Vert\td_{X_1,\cdots 
X_N}^{g_1}R^{g_1}-\td_{X_1,\cdots X_N}^{g_2}R^{g_2}\Vert_{g_1}$ is bounded
for any $g_1, g_2$ chosen among $g_B,$ $g_{KE}$ and $g_{M}.$  Hence $g_B,
g_{KE}$ and $g_M$ are all comparable up to any order of derivatives.

\bs
\begin{center}
{\bf  References}
\end{center}

\ms
\ni [At] Atiyah, M. F., Elliptic operators, discrete groups, and von Neumann
algebras, Ast\'erisque 32-33(1976), 43-72.

\ms
\ni [Au]
 Aubin, T., Some nonlinear problems in Riemannian geometry. Springer Monographs in Mathematics. Springer-Verlag, Berlin, 1998.
  
\ms
\ni
[B] Borel, A., 
Compact Clifford-Klein forms of symmetric spaces, 
Topology 2 (1963), 111-122.

\ms
\ni
[CG] Cheeger, J., Gromov, M., On the characteristic numbers of complete
manifolds of bounded curvature and finite volume, Differential geometry
and complex analysis, 115-154, Springer, Berlin, 1985.

\ms
\ni[De] Demailly, J.-P., Champs magnetiques et inegalites de Morse pour la $d"$-cohomologie.  35 (1985), no. 4, 189--229.

\ms
\ni
[DF] Diederich, K., Fornaess, J. E., Pseudoconvex domains: bounded strictly plurisubharmonic exhaustion functions. Invent. Math. 39 (1977), no. 2, 129--141.

\ms
\ni [Do] Donnelly, H., $L^2$ cohomology of the Bergman metric for weakly pseudoconvex domains.  Illinois J. Math.  41  (1997),  no. 1, 151--160.

\ms
\ni
[Ga] Gardiner, F., Teichm\"uller theory and quadratic differentials, 
Wiley Interscience, New York, 1987.

\ms
\ni
[Gr] Gromov, M., K\"ahler hyperbolicity and $L^2$-Hodge theory.  J. Differential Geom.  33  (1991),  no. 1, 263--292.

\ms
\ni[HZ] Harer, J., Zagier, D.  The Euler characteristic of the moduli space of curves. Invent. Math. 85 (1986), no. 3, 457--485.

\ms
\ni
[He] Helgason, S., Differential geometry, Lie groups and symmetric spaces, 
Academic Press 1978
 
\ms
\ni
[KM] Knudsen, F. and Mumford, D., The projectivity of the moduli space of stable curves. I. 
Preliminaries on "det" and "Div".  Math. Scand.  39  (1976), no. 1, 19--55. 

\ms
\ni
[Ko] Koll\'ar, J., Shafarevich maps and automorphic forms. 
M. B. Porter Lectures. Princeton University Press, Princeton, NJ, 1995. 
 
\ms
\ni [Kr] Krantz, S., Function theory of several complex variables. Second edition. The Wadsworth \& Brooks/Cole Mathematics Series. Wadsworth \& Brooks/Cole Advanced Books \& Software, Pacific Grove, CA, 1992. xvi+557 pp

\ms
\ni
[Mc] McMullen, C. T., The moduli space of Riemann surfaces is K\"ahler hyperbolic, Ann. of Math.
151 (2000), 327--357.

\ms
\ni
[Mo] Mok, N., The Serre problem on Riemann surfaces. Math. Ann. 258 (1981/82), no. 2, 145--168. 

\ms
\ni
[MY] Mok, N.,  Yau, S., Completeness of the K\"ahler-Einstein metric on bounded domains and the characterization of domains of holomorphy by curvature conditions.  The mathematical heritage of Henri Poincar\'e, Part 1 (Bloomington, Ind., 1980),  41--59, Proc. Sympos. Pure Math., 39, Amer. Math. Soc., Providence, RI, 1983.

\ms
\ni
[MZ] Mok, N., Zhong, J. Q., Compactifying complete Kahler-Einstein manifolds of finite topological type and bounded curvature. Ann. of Math. (2) 129 (1989), no. 3, 427--470.
 
\ms
\ni 
[NT] Nadel, A., Tsuji, H., Compactification of complete Kahler manifolds of negative Ricci curvature. J. Diff. Geom. 28 (1988), no. 3, 503--512. 
 
\ms
\ni
[Ro] Royden, H., The Ahlfors-Schwarz lemma in several complex variables,
 Comment. Math. Helv.  55  (1980), no. 4, 547--558. 

\ms
\ni
[Se] Serre, J.-P., Quelques problèmes globaux relatifs aux vari\'et\'es de Stein, Colloque sur les fonctions de plusieurs variables, tenu \'a Bruxelles, 1953, pp. 57--68. Georges Thone, Li\`ege; Masson \& Cie, Paris, 1953.

\ms
\ni
[Si] Siu, Y.-T., Holomorphic fiber bundles whose fibers are bounded Stein domains with zero first Betti number. Math. Ann. 219 (1976), no. 2, 171--192.

\ms
\ni
[Sk] Skoda, H., Fibr\'es holomorphes à base et \`a fibre de Stein,
Invent. Math. 43 (1977), no. 2, 97--107.

\ms
\ni
[St]  Stehl\'e, J.-L., Fonctions plurisousharmoniques et convexit\'e holomorphe de certains fibr\'es analytiques. (French) C. R. Acad. Sci. Paris S\'er. A 279 (1974), 235--238. 

\ms
\ni
[W] Wolpert, S., Weil-Petersson perspectives, arXiv: math.DG/0502519.

\ms
\ni
[Y1] Yeung, S.-K., Compactification of Kahler manifolds with negative Ricci curvature. Invent. Math. 106 (1991), no. 1, 13--25. 

\ms
\ni
[Y2] Yeung, S.-K., Bounded smooth strictly plurisubharmonic exhaustion functions 
on Teichm\"uller spaces, Math. Res. Letters, (10) 2003, 391-400.

\ms
\ni
[Y3] Yeung, S.-K., Quasi-isometry of metrics on Teichm\"uller spaces, 
Int. Math. Res. Not. 2005, 239-255. 

\ms
\ni
[Y4] Yeung, S.-K., Bergman metric on Teichm\"uller spaces and moduli spaces of curves, in Recent progress on some problems in several complex variables
and partial differential equations,  Contempory Mathematics 400. 

\end{document}